\documentclass[10pt]{article}
\usepackage{subeqn,amsfonts}
\newcommand{\sect}[1]{\setcounter{equation}{0}\section{#1}}

\newcommand{\be}{\begin{equation}}
\newcommand{\ee}{\end{equation}}
\newcommand{\bea}{\begin{eqnarray}}
\newcommand{\eea}{\end{eqnarray}}

\newcounter{remark}

\newcommand{\cA}{{\cal A}}

\newcommand{\cR}{{\cal R}}
\newcommand{\cY}{{\cal Y}}

\newcommand{\II}{{\mathbb I}}
\newcommand{\ZZ}{{\mathbb Z}}
\newcommand{\eps}{{\varepsilon}}
\newcommand{\elp}{{\cA_{q,p}(\widehat{sl}(2)_c)}}
\newcommand{\snh}{{\mbox{snh}}}

\newcommand{\qpvir}{{{\cal W}_{q,p}(sl(2))}}
\newcommand{\remark}{{\sl Remark \addtocounter{remark}{1}\theremark}}
\newtheorem{defin}{Definition}
\newtheorem{prop}{Proposition}
\newtheorem{thm}{Theorem}

\begin{document}
\newpage
\pagestyle{plain}
\setcounter{remark}{0}
\vfill
\begin{center}

{{\bf From quantum to elliptic algebras}}

\vspace{7mm}

{ J. Avan}

\vspace{4mm}

{\it LPTHE, CNRS-URA 280, Universit\'es Paris VI/VII, France}

\vspace{7mm}

{ L. Frappat \footnote{On leave of absence from Laboratoire de
Physique Th\'eorique  ENSLAPP.}}

\vspace{4mm}

{\it Centre de Recherches Math\'ematiques, Universit\'e de Montr\'eal,
Canada}

\vspace{7mm}

{ M. Rossi, P. Sorba}

\vspace{4mm}

{\it Laboratoire de Physique Th\'eorique ENSLAPP \footnote{URA 1436 du
CNRS, associ\'ee \`a l'\'Ecole Normale Sup\'erieure de Lyon et \`a
l'Universit\'e de Savoie.}\\
Annecy-le-Vieux C\'edex
and ENS Lyon, Lyon C\'edex 07, France}

\end{center}

\vfill

\begin{abstract}
It is shown that the elliptic algebra $\elp$ at the critical level $c=-2$
has a multidimensional center containing some trace-like operators
$t(z)$.
A family of Poisson structures indexed by a non-negative integer and
containing the $q$-deformed Virasoro algebra is constructed on this
center.
We show also that $t(z)$ close an exchange algebra when $p^m=q^{c+2}$ for
$m\in\ZZ$, they commute when in addition
$p=q^{2k}$ for $k$ integer non-zero, and they belong to the center of
$\elp$ when $k$ is odd. The Poisson structures obtained for $t(z)$ in
these classical limits contain the $q$-deformed Virasoro algebra,
characterizing the structures at $p \neq q^{2k}$ as new $\qpvir$
algebras.
\end{abstract}
\vspace{1truecm}
{\it Talk presented by M.Rossi at the 6th International Colloquium
``Quantum Groups and Integrable Systems'', Prague, 19-21 June 1997.}

\vfill
\vfill

\rightline{ENSLAPP-AL-657/97}
\rightline{CRM-2485}
\rightline{PAR-LPTHE 97-37}
\rightline{q-alg/9707034}
\rightline{July 1997}

\newpage
\pagestyle{plain}

\sect{Introduction}

The concept of $q$-deformed Virasoro and $W_N$ algebras has recently
arisen in
connection with some aspects of integrable systems. In particular these
algebras have been introduced \cite{SKAO,AKOS} as an extension of
the Virasoro and $W_N$ algebras identified in the quantum Calogero--Moser
model \cite{AMOS}. In the same way as Jack polynomials (eigenfunctions of
the quantum Calogero--Moser model) arise as singular vectors of $W_N$
algebras \cite{AMOS,MaChe}, MacDonald polynomials (eigenfunctions of
the trigonometric Ruijsenaars--Schneider model) are singular vectors of
the $q$-deformed $W_N$ algebras.

These deformed algebras were shown to arise in fact from a procedure of
construction mimicking the already known scheme \cite{FF1} for undeformed
Virasoro and $W_N$ algebras: at the critical value $c=-N$, the quantum
affine Lie algebra $U_q(\widehat{sl}(N)_c)$ has a
multidimensional center \cite{RSTS} where a $q$-deformed Poisson bracket
can be
defined as limit of the commutator structure \cite{FR}. This Poisson
bracket is the semi-classical limit of the $q$-deformed Virasoro (for
$N=2$) or $W_N$ algebras. Its quantization, performed in \cite{FF2}, 
can be considered as the most
general definition of $q$-deformed $W_N$ algebras, often called
${{\cal W}_{q,p}(sl(N))}$, $p$ being related to the quantization
parameter. In fact, the construction of classical and quantized algebras
was not achieved in \cite{FR,FF2} by direct computation but using the
$q$-deformed bozonization \cite{AOS} of $U_q(\widehat{sl}(N)_c)$.
Interestingly the $q$-deformed $W_N$ algebras are characterized by
elliptic structure coefficients: for instance, the $q$-deformed Virasoro
algebra of \cite{SKAO} is defined by the generating operator $T(z)$ such
that
\be
f_{1,2}(w/z) \, T(z) \, T(w) - f_{1,2}(z/w) \, T(w) \, T(z)
= \frac{(1-q)(1-p/q)}{1-p} \left( \delta\Big(\frac{w}{zp}\Big)
- \delta\Big(\frac{wp}{z}\Big) \right)\, ,
\ee
where
\be
f_{1,2}(x) = \frac{1}{1-x} ~ \frac{(x|q,pq^{-1};p^2)_\infty}
{(x|pq,p^2q^{-1};p^2)_\infty} \, , \quad
(x|a_1 \dots a_k;t)_\infty \equiv \prod_{i=1}^k \prod_{n=0}^\infty
(1-a_ixt^n)
\, .
\ee

The parameters $p$ and $q$ are rewritten as $q=e^h$, $p=e^{h(1-\beta)}$:
$h$ is the deformation parameter and $\beta$ is the ``quantization''
parameter, the semi-classical limit $\beta \rightarrow 0$ giving back the
$q$-deformed Poisson bracket of \cite{FR} and the limit $h \rightarrow 0$
giving back the
linear Virasoro algebra.

Two natural problems arise in this context. The first is the extension of
these Poisson brackets constructions to the elliptic quantum algebra
$\elp$ \cite{FIJKMY,KLP,JMK}, which is a double deformation of $U(\widehat
{sl}(2)_c)$, the limit $p\rightarrow 0$ giving the quantum affine algebra
$U_q(\widehat {sl}(2)_c)$. In this context we will show the existence of a
multidimensional center at $c=-2$, where we will construct a set of
Poisson brackets containing the $q$-deformed Virasoro algebra introduced
by \cite{FR}.

A second question is the analysis of the connection between $\elp$ and
quantized $q$-deformed Virasoro algebra, both of which depend on elliptic
structure functions. Our result is that if $q^{c+2}=p^m$ for any $m \in
\ZZ \backslash \{0\}$, the algebra $\elp$ contains a quadratic subalgebra
which builds a
natural quantization of the Poisson bracket structure of the $q$-deformed
Virasoro algebra; in such a way we will construct a family of $\qpvir$
algebras in the framework of the elliptic quantum algebra $\elp$.

\sect{The elliptic quantum algebra $\elp$}

Consider the $R$-matrix of the eight vertex model found by Baxter
\cite{Ba}:
\be
R_{12}(x) = \frac{1}{\mu(x)} \left(\begin{array}{cccc}
a(u) & 0 & 0 & d(u) \cr 0 & b(u) & c(u) & 0 \cr
0 & c(u) & b(u) & 0 \cr d(u) & 0 & 0 & a(u) \cr
\end{array}\right) \label{eq21}
\ee
where the functions $a(u), b(u), c(u), d(u)$ are given by
\be
a(u) = \frac{\snh(\lambda-u)}{\snh(\lambda)} \,, \quad
b(u) = \frac{\snh(u)}{\snh(\lambda)} \,, \quad
c(u) = 1 \,, \quad
d(u) = k\,\snh(\lambda-u)\snh(u) \,.\label{eq22}
\ee
The function $\snh(u)$ is defined by $\snh(u) = -i\mbox{sn}(iu)$ where
$\mbox{sn}(u)$ is Jacobi's elliptic function with modulus $k$.
The functions (\ref{eq22}) can be seen as depending on:
\be
p = \exp\Big(-\frac{\pi K'}{K}\Big) \,, \qquad
q = - \exp\Big(-\frac{\pi\lambda}{2K}\Big) \,, \qquad
x = \exp\Big(\frac{\pi u}{2K}\Big) \,,
\ee
where $K$, $K^\prime$ denote the elliptic integrals \cite{Ba}.
The factor $\mu(x)$ is given by \cite{JMK}:
\bea
&& \frac{1}{\mu(x)} = \frac{1}{\kappa(x^2)} \frac{(p^2;p^2)_\infty}
{(p;p)_\infty^2} \frac{\Theta_{p^2}(px^2)\Theta_{p^2}(q^2)}
{\Theta_{p^2}(q^2x^2)} \,, \\
&& \frac{1}{\kappa(x^2)} = \frac{(q^4x^{-2};p,q^4)_\infty
(q^2x^2;p,q^4)_\infty (px^{-2};p,q^4)_\infty (pq^2x^2;p,q^4)_\infty}
{(q^4x^2;p,q^4)_\infty (q^2x^{-2};p,q^4)_\infty (px^2;p,q^4)_\infty
(pq^2x^{-2};p,q^4)_\infty} \,,
\eea
where
\bea
(x;p_1,\dots,p_m)_\infty = \prod_{n_i \ge 0} (1-xp_1^{n_1} \dots
p_m^{n_m}) \\
\Theta_{p^2}(x) = (x;p^2)_\infty (p^2x^{-1};p^2)_\infty (p^2;p^2)_\infty
\,.
\eea
To avoid singularities in the functions (2.2,4) we will suppose $|q|, \,
|p|<1$.
\begin{prop}\label{propone}
The matrix $R_{12}$ has the following properties:
\\
\indent -- unitarity: $R_{21}(x^{-1}) R_{12}(x) = 1$,
\\
\indent -- crossing symmetry: $R_{21}(x^{-1})^{t_1}
= (\sigma^1 \otimes \II) R_{12}(-q^{-1}x) (\sigma^1 \otimes \II)$,
\\
\indent -- antisymmetry: $R_{12}(-x) = - (\sigma^3 \otimes \II)
R_{12}(x) (\sigma^3 \otimes \II)$,
\\
where $\sigma^1,\sigma^2,\sigma^3$ are the Pauli matrices and
$t_i$ is the transposition in the space $i$.
\end{prop}
The proof is straightforward by direct calculation.

For the definition of the elliptic quantum algebra $\elp$, we need to
use a modified $R$-matrix $R_{12}^+(x)$ defined by
\be
R_{12}^+(x) = \tau(q^{1/2}x^{-1}) R_{12}(x) \, , \quad
\tau(x) = x^{-1} \frac{\Theta _{q^4}(qx^2)}{\Theta _{q^4}(qx^{-2})} \,.
\ee
The matrix $R_{12}^+(x)$ obeys a quasi-periodicity property:
\be
R_{12}^+(-p^{\frac {1}{2}}x)= (\sigma^1 \otimes \II)
\left (R_{21}^+(x^{-1})\right )^{-1} (\sigma^1 \otimes \II) \, .
\ee

\medskip

The elliptic quantum algebra $\elp$ has been introduced by \cite{FIJKMY}.
It is an algebra of operators $L_{\eps\eps',n}^\pm$ ($\eps,\eps'=+$ or
$-$) such that $L_{\eps\eps',n}^\pm = 0 \mbox{ if } \eps\eps' \ne (-1)^n$
and, defining $L_{\eps\eps'}^\pm(z) = \sum_{n\in\ZZ}
L_{\eps\eps',n}^\pm z^n$ and encapsulating them into $2 \times 2$ matrices
$L^{\pm}(z)$, one has, with the definition $R^{+*}_{12}(x,q,p)=
R^+_{12}(x,q,pq^{-2c})$:
\be
\begin{array}{l}
R^+_{12}(z/w) ~ L^\pm_1(z) ~ L^\pm_2(w) =
L^\pm_2(w) ~ L^\pm_1(z) ~ R^{+*}_{12}(z/w) \,, \\
\bigg. R^+_{12}(q^{c/2}z/w) ~ L^+_1(z) ~ L^-_2(w) =
L^-_2(w) ~ L^+_1(z) ~ R^{+*}_{12}(q^{-c/2}z/w) \,, \\
q-{\mbox {det}}L^+(z)\equiv
L^+_{++}(q^{-1}z)L^+_{--}(z)-L^+_{-+}(q^{-1}z)L^+_{+-}(z)=q^{\frac {c}{2}}
\, , \\
L_{\eps\eps'}^-(z)=\eps \eps' L_{-\eps, -\eps'}^+(p^{\frac
{1}{2}}q^{-\frac {c}{2}}z) \, .
\end{array} \label{eq216}
\ee
We now state the first important result of our study.

\sect{The center of $\elp$}

\begin{thm}\label{thmone}
At $c=-2$, the operators generated by
\be
t(z) = {\mbox Tr}(L(z)) = {\mbox Tr}\Big(L^+(q^{c/2}z) L^-(z)^{-1}\Big)
\label{eq31}
\ee
commute with the algebra $\elp$ and then belong to its center.
\end{thm}
The formula for $t(z)$ in the elliptic case is exactly identical to the
one in
the trigonometric case \cite{RSTS}. The proof follows on similar lines
using explicitly the crossing symmetry and the unitarity of the
$R$-matrix.

We now study the specific behaviour of the exchange algebra of $t(z)$ with
$t(w)$ in the neighborhood of $c=-2$.

\sect{Poisson algebra of $t(z)$}

By virtue of Theorem \ref{thmone} the elements $t(z)$ and $t(w)$ are
mutually commuting at the
critical level $c=-2$. This implies a natural Poisson structure on the
algebra generated by them: if $\Big[ t(z),t(w) \Big] = (c+2) \ell(z,w) +
o(c+2)$, then a Poisson bracket is yielded by
$\Big\{ t(z)_{cr},t(w)_{cr} \Big\} = \ell(z,w)_{cr}$
(``cr'' means that all expressions are taken at $c=-2$).
We have the following result:
\begin{thm}\label{thmtwo}
Under its natural Poisson bracket at $c=-2$ the algebra generated by
$t(z)$ is closed. One has indeed (we suppress the subscript ``$cr$'' and
define $x=z/w$):
\be
\Big\{ t(z),t(w) \Big\} = -(\ln q) \left(x^{-1}\frac{d}{dx^{-1}}
\ln\tau(q^{1/2}x^{-1}) - x\frac{d}{dx}\ln\tau(q^{1/2}x)\right)
~ t(z)t(w) \,. \label{eq41}
\ee
\end{thm}
{\bf Proof:}
From the definition of the element $t(z)$, one has
\be
t(z)t(w) = L(z)^{i_1}_{i_1} \, L(w)^{i_2}_{i_2} =
L^+(q^{\frac {c}{2}}z)^{j_1}_{i_1} (L^-(z)^{-1})^{i_1}_{j_1} \,
L^+(q^{\frac {c}{2}}w)^{j_2}_{i_2}  (L^-(w)^{-1})^{i_2}_{j_2} \,.
\ee
The exchange relations of (\ref{eq216})
and the properties of the matrix $R_{12}(x)$ given in proposition
\ref{propone} allow us to move the matrices $L^+(q^{\frac {c}{2}}w)$,
$L^-(w)^{-1}$ to the left of the matrices $L^+(q^{\frac {c}{2}}z)$,
$L^-(z)^{-1}$. One obtains
\be
t(z)t(w) = \cY(z/w)^{i_1i_2}_{j_1j_2} ~ L(w)^{j_2}_{i_2} ~
L(z)^{j_1}_{i_1} \,,
\label{eq44}
\ee
where the matrix $\cY(z/w)$ is factorized in the following way:
\be
\cY(z/w) = T(z/w) \cR(z/w) \,. \label{eq46}
\ee
The matrix factor $\cR(z/w)$ depends only on the matrix (\ref{eq21}):
\be
\cR(z/w) = \bigg(\Big(R_{12}(w/z)\, R_{12}(q^{-c-2}z/w)
\, R_{12}(w/z)\Big)^{t_2} \, {R_{12}(q^cz/w)}^{t_2}\bigg)^{t_2} \, ,
\ee
while the numerical prefactor $T(z/w)$ contains only a $\tau$ dependence:
\be
T(z/w) = \frac{\tau(q^{1/2}z/w)\tau(q^{-c+1/2}w/z)}
{\tau(q^{-c-3/2}z/w)\tau(q^{1/2}w/z)} \,. \label{eq48}
\ee
One easily checks the nice behaviour of $T(z/w)$ and $\cR(z/w)$
at $c=-2$:
\be
T(z/w)_{cr} = 1  \, , \quad  \cR(z/w)_{cr} = \II_2 \otimes \II_2
\qquad \Longrightarrow \qquad
\cY(z/w)_{cr} = \II_2 \otimes \II_2 \,. \label{eq49}
\ee
\medskip
One then computes the Poisson structure from the exchange algebra
(\ref{eq44})
in the neighborhood of $c=-2$. From equations (\ref{eq44}) and
(\ref{eq49})
one writes
\be
t(z)t(w) = t(w)t(z) +
(c+2)\left(\frac{d\cY}{dc}(z/w)\right)^{i_1i_2}_{j_1j_2}
~ L(w)^{j_2}_{i_2} ~ L(z)^{j_1}_{i_1}\bigg\vert_{cr} + o(c+2)
\ee
and therefore
\be
\Big\{ t(z),t(w) \Big\} =
\left(\frac{d\cY}{dc}(z/w)\right)^{i_1i_2}_{j_1j_2}
~ L(w)^{j_2}_{i_2} ~ L(z)^{j_1}_{i_1}\bigg\vert_{cr} \, . \label{eq411}
\ee
The equations (\ref{eq46}) and (\ref{eq49}) imply
\be
\frac{d\cY}{dc}(x)\bigg\vert_{cr} = \frac{dT}{dc}(x)\bigg\vert_{cr}
\II_2 \otimes \II_2 ~+~ \frac{d\cR}{dc}(x)\bigg\vert_{cr} \,.
\label{eq412}
\ee
After a long calculation, using various tricks in elliptic
functions theory, one shows that:
\be
\frac{d \cR}{dc}(x)\bigg\vert_{cr} =0 \,,
\label{eq413}
\ee
\be
\frac{dT}{dc}(x)\bigg\vert_{cr} = -(\ln q) \left(x^{-1}\frac{d}{dx^{-1}}
\ln\tau(q^{1/2}x^{-1}) - x\frac{d}{dx}\ln\tau(q^{1/2}x)\right)
\,.\label{eq414}
\ee
From equations (\ref{eq411}-\ref{eq414}) formula (\ref{eq41}) of
Theorem \ref{thmtwo} immediately follows. \hfill \rule{5pt}{5pt}

\medskip

The structure of the Poisson bracket (\ref{eq41}) derives wholly from the
$\tau$ factor: so any dependence in $p$ is absent in its structure
function.

\medskip

{From} the equation (\ref{eq41}) and the definition of
$\tau(x)$, one gets easily
\bea
\Big\{ t(z),t(w) \Big\} = -(2\ln q) \left[ \sum_{n \ge 0} ~ \left(
\frac{2x^2q^{4n+2}}{1-x^2q^{4n+2}} -
\frac{2x^{-2}q^{4n+2}}{1-x^{-2}q^{4n+2}}
\right ) \right. \nonumber \\
 \left. +\sum_{n>0} ~ \left(
- \frac{2x^2q^{4n}}{1-x^2q^{4n}} +
\frac{2x^{-2}q^{4n}}{1-x^{-2}q^{4n}}\right)
- \frac{x^2}{1-x^2} + \frac{x^{-2}}{1-x^{-2}} \right] ~
t(z) \, t(w) \, , \label{eq415}
\eea
where $x=z/w$.
Interpretation of the formula (\ref{eq415}) must now be given in terms of
the modes $t_n$ of $t(z)$, defined by:
\be
t_n = \oint_C \frac{dz}{2\pi iz} \, z^{-n} \, t(z) \, .
\ee
The structure function $f(z/w)$ which defines the Poisson bracket
(\ref{eq415}) is periodic
with period $q^2$ and has simple poles at $z/w = \pm q^k$ for $k \in
\ZZ$.
In particular it is singular at $z/w = \pm 1$.
As a consequence, the expected definition of the Poisson bracket
$\{t_n,t_m\}$ as a double contour integral of (\ref{eq415}) must be made
more precise. Deformation of, say, the $w$-contour while the $z$-contour
is
kept fixed may induce the crossing of singularities of $f(z/w)$ which in
turn
modifies the computed value of the Poisson bracket. In particular
the singularity at $z/w= \pm 1$ implies that one cannot identify
a double contour integral with its permuted. As a consequence the quantity
$\displaystyle{\oint_{C_1} \frac{dz}{2\pi iz} \oint_{C_2}
\frac{dw}{2\pi iw} \, z^{-n} \, w^{-m} \, f(z/w) \, t(z) \, t(w)}$
is {\em not} antisymmetric under the exchange $n \leftrightarrow m$
and cannot be taken as a Poisson bracket.
This leads us to define the Poisson bracket as:
\begin{defin}\label{defone}
\be
\{ t_n,t_m \} = \frac{1}{2} \left( \oint_{C_1} \frac{dz}{2\pi iz}
\oint_{C_2}
\frac{dw}{2\pi iw} + \oint_{C_2} \frac{dz}{2\pi iz} \oint_{C_1}
\frac{dw}{2\pi iw} \right) z^{-n} \, w^{-m} \, f(z/w) \, t(z) \, t(w) \,
.
\label{eq416}
\ee
\end{defin}
Such a procedure guarantees the antisymmetry of the postulated Poisson
structure due to the property $f(z/w) = -f(w/z)$.

The presence of singularities at
$z/w = \pm q^k$ where $k \ne 0$ introduce a dependence of the Poisson
bracket (\ref{eq416}) on the domains of integration. If we choose the
contours $C_1$ and
$C_2$ to be circles of radii $R_1$ and $R_2$ respectively, the following
proposition holds:
\begin{prop}\label{proptwo}
For any $k \in \ZZ^+$ such that $R_1/R_2 \in [q^{\pm k},q^{\pm(k+1)}]$,
Definition \ref{defone} defines a consistent (that is ``antisymmetric and
obeying the Jacobi identity'') Poisson bracket whose specific
form, depending on $k$, is:
\end{prop}
\bea
&& \{ t_n,t_m \}_k = (-1)^{k+1} 2 \ln q \oint_{C_1} ~ \frac{dz}{2\pi iz}
\oint_{C_2}~ \frac{dw}{2\pi iw} \cdot \nonumber \\
&& \cdot \sum_{s \in \ZZ} \frac{q^{(2k+1)s} - q^{-(2k+1)s}}
{q^s+q^{-s}} \Big( \frac{z}{w} \Big)^{2s} z^{-n} \, w^{-m} \, t(w) \, t(z)
\, .
\label{eq417}
\eea
We observe that the form of the Poisson brackets (\ref{eq417}) is similar
to the form of the Poisson bracket obtained by \cite{FR}. In particular
our Poisson bracket at $k=1$ is the one in \cite{FR} where the purely
central term is multiplied by
$t(z)t(zq)$. However one has to remember that in \cite{FR} a
particular representation of $U_q(\widehat{sl}(2)_c)$ in terms of
quasi-bosons is used. It is possible that an analogous bosonization of
$\elp$ leads to a degeneracy of such terms as $t(z)t(zq)$ giving the
result of \cite{FR}:
unluckly at this time a bosonized version of $\elp$ is
available only at $c=1$ \cite{FSHY}, using bosonized vertex operators
constructed in \cite{AJMP}.
\medskip

\sect{Quadratic subalgebras in $\elp$}

We now turn to the task of identifying possible connections between
$\elp$ and $\qpvir$. We first prove:
\begin{thm}
If $p$, $q$, $c$ are connected by the relation $p^m=q^{c+2}$, $m \in \ZZ$,
the operators $t(z)$ realize an exchange algebra with all generators
$L^\pm (w)$ of $\elp$:
\be
t(z)L^+(w) = F\Big(m,q^{\frac {c}{2}}\frac{w}{z}\Big) L^+(w)t(z) \,,
\quad
t(z)L^-(w) = F\Big(m,-p^{\frac {1}{2}}\frac{w}{z}\Big) L^-(w)t(z) \,,
\label{eq51}
\ee
where
\begin{subequations}
\label{eq52}
\bea
F(m,x) &=& \prod_{s=1}^{2m} q^{-1} \,
\frac{\Theta_{q^4}(x^{2}q^2p^{-s}) \, \Theta_{q^4}(x^{-2}q^2p^s)}
{\Theta_{q^4}(x^{-2}p^s) \, \Theta_{q^4}(x^2p^{-s})}
\quad \mbox{for $m>0$} \,, \\
F(m,x) &=& \prod_{s=0}^{2|m|-1} q \,
\frac{\Theta_{q^4}(x^2p^s) \, \Theta_{q^4}(x^{-2}p^{-s})}
{\Theta_{q^4}(x^2q^2p^s) \, \Theta_{q^4}(x^{-2}q^2p^{-s})}
\quad \mbox{for $m<0$} \,.
\eea
\end{subequations}
\end{thm}

The proof is easy to perform using the definition of $\elp$ and the
properties
(especially the quasi-periodicity) of $R$.

\remark : For $m=0$, the relation can be realized in two ways: either
$c=-2$, which is the case studied in chapter 3 and leads directly to a
central $t(z)$ ($F(m,x)=1$); or $q=\exp \left
({\frac{2i\pi\ZZ}{c+2}}\right )$,
hence $|q|=1$, which we have decided not to consider owing to the
singularities in the elliptic functions defining $\elp$. Hence $m=0$ will
be
disregarded from now on.

\medskip

An immediate corollary is:
\begin{thm}
When $p^m=q^{c+2}$, $t(z)$ closes a quadratic subalgebra:
\be
t(z)t(w) = {\cal Y}_{p,q,m}\Big(\frac{w}{z}\Big) \, t(w)t(z) \label{eq53}
\ee
where
\be
{\cal Y}_{p,q,m}(x) = \left\{ \begin{array}{ll}
\displaystyle \left[ \prod_{s=1}^{2m-1} x^{2}
\frac{\Theta_{q^4}(x^{-2}p^s) \, \Theta_{q^4}(x^2q^2p^s)}
{\Theta_{q^4}(x^2p^s) \, \Theta_{q^4}(x^{-2}q^2p^s)} \right]^2
& \mbox{for $m>0$} \,, \\ \\
\displaystyle \left[ \prod_{s=1}^{2|m|} x^{2}
\frac{\Theta_{q^4}(x^{-2}p^s) \, \Theta_{q^4}(x^2q^2p^s)}
{\Theta_{q^4}(x^2p^s) \, \Theta_{q^4}(x^{-2}q^2p^s)} \right]^2
& \mbox{for $m<0$} \,. \\
\end{array} \right.
\label{eq54}
\ee
\end{thm}

\medskip

\noindent
The proof is obvious from (\ref {eq51}) and the definition (2.7) of
$\Theta _a(x)$.

\remark: When $m=1$ the exchange function in (\ref{eq53}) is exactly
the square of the exchange function in the quantized $q$-deformed
Virasoro algebra proposed in \cite{SKAO}, once the replacements $q^2
\rightarrow p$, $p \rightarrow q$, $x^2 \rightarrow x$ are done.
This is a first indication in our context that the elliptic algebra $\elp$
appears to be
the natural setting to define quantized $q$-deformed Virasoro algebras.

\remark: As an additional connection we notice that all
functions ${\cal Y}_{p,q,m}(x)$ obey the Feigin-Frenkel identities
\cite{FF2} for the exchange function of \cite{SKAO}:
\be
{\cal Y}(xq^2) = {\cal Y}(x) \,, \quad {\cal Y}(xq) = {\cal Y}(x^{-1})
\,.
\label{eq55}
\ee
Our exchange algebras then appear as natural candidates for $\qpvir$
algebras generalizing the one of \cite{SKAO}. This interpretation will be
reinforced by the next results.
First of all we state the following theorem:

\begin{thm} \label{thm6}
For $p=q^{2k}$, $k\in \ZZ \backslash \{ 0 \}$, one has
\bea
F(m,x) &=& 1 \qquad \mbox{for $k$ odd} \,, \\
F(m,x) &=& q^{-2m}x^{4m}\left[ \frac{\Theta _{q^4}(x^2q^2)}{\Theta
_{q^4}(x^2)}
\right]^{4m} \qquad \mbox{for $k$ even} \,.
\eea
Hence when $k$ is odd $t(z)$ is in the center of the algebra $\elp$, while
when $k$ is even $t(z)$ is {\sl not} in a (hypothetical) center of
$\elp$.
However {\sl in both cases}, one has $[t(z) \, , \, t(w)] =0$.
\end{thm}

\noindent
{\bf Proof:}
Theorem \ref{thm6} is easily proved using the explicit expression for
$F(m,x)$ and the definition (2.7) of $\Theta$-functions.
The case $k=0$ is excluded since it would lead to $p=1$ and singularities
in the definition of $\elp$.
\hfill \rule{5pt}{5pt}

\medskip

This now allows us to define Poisson structures even though $t(z)$ is not
in the center of $\elp$ for $k$ even. They
are obtained as limits of the exchange algebras (\ref{eq53}). Since the
initial non-abelian structure for $t(z)$ is closed, the exchange algebras
(\ref{eq53}) are natural quantizations of the Poisson algebras which we
obtain.

\begin{thm}
Setting $q^{2k}=p^{1-\frac{\beta }{2}}$ for any integer $k \ne 0$, one
defines in the limit $\beta \rightarrow 0$ the following Poisson
structures ($x=z/w$):
\begin{subequations}
\label{eq56}
\bea
&&\Big\{ t(z) \, , \, t(w)\Big\} \equiv \lim_{\beta \rightarrow 0} \,
\frac{1}{\beta} \, \left[ t(z)t(w)-t(w)t(z) \right] \nonumber \\
&&= 2km \, \ln q \left\{
-\frac{x^2}{1-x^2} +\frac{x^{-2}}{1-x^{-2}} + \sum_{n=0}^{\infty}
\left [ \frac{2x^2q^{4n}}{1-x^2q^{4n}} -
\frac{2x^2q^{4n+2}}{1-x^2q^{4n+2}}
\right. \right . \nonumber \\
&& \hspace{20mm} \left. \left. - \frac{2x^{-2}q^{4n}}{1-x^{-2}q^{4n}}
+\frac{2x^{-2}q^{4n+2}}{1-x^{-2}q^{4n+2}} \right ] \right\} t(z)t(w)
\quad \mbox{for $k$ odd} \,, \\
&&= -2km(2m-1)\, \ln q \left\{
-\frac{x^2}{1-x^2} + \frac{x^{-2}}{1-x^{-2}} + \sum_{n=0}^{\infty}
\left [ \frac{2x^2q^{4n}}{1-x^2q^{4n}} -
\frac{2x^2q^{4n+2}}{1-x^2q^{4n+2}}
\right. \right. \nonumber \\
&& \hspace{20mm} \left. \left. - \frac{2x^{-2}q^{4n}}{1-x^{-2}q^{4n}}
+ \frac{2x^{-2}q^{4n+2}}{1-x^{-2}q^{4n+2}} \right ] \right\}t(z)t(w)
\quad \mbox{for $k$ even} \,.
\eea
\end{subequations}
\end{thm}

\noindent
{\bf Proof:}
We note that
\be
\Big\{ t(z) \, , \, t(w) \Big\} = \frac{d{\cal Y}_{p,q,m}}{d\beta}
\left (\frac {w}{z}\right)\Big\vert_{\beta=0} \, t(z)t(w) =
\frac{d\ln{\cal Y}_{p,q,m}}{d\beta}\left (\frac {w}{z}\right)
\Big\vert_{\beta=0} \, t(z)t(w) \,,
\ee
the two equalities coming from the fact that ${\cal Y}_{p,q,m} = 1$ when
$q^{2k}=p$.
The proof is then obvious from (5.3,4) and the definition of
$\Theta$-functions as absolutely convergent products (for $|q|<1$),
hence the series in (\ref{eq56}) are convergent and define
univocally a structure function for $\{ t(z)\, , \, t(w) \}$.
\hfill \rule{5pt}{5pt}

\medskip

The formula (\ref{eq56}) coincides exactly with the Poisson structure of
$t(z)$ at $c=-2$, provided one reabsorbs $km$ and $-km(2m-1)$ into the
definition of the classical limit as $\beta \rightarrow km\beta$ for $k$
odd and $\beta \rightarrow -km(2m-1)\beta$ for $k$ even.
So Theorem 6 provides us with an immediate interpretation of the
quadratic structures (\ref{eq53}). Since we have seen that
the Poisson structures derived from (\ref{eq415}) contained in particular
the
$q$-deformed Virasoro algebra (up to the delicate point of the central
extension which is not explicit in (\ref{eq415})), the quadratic algebras
(\ref{eq53}) are inequivalent (for different values of $m$ !)
quantizations
of the classical $q$-deformed Virasoro algebra, globally defined on the
$\ZZ$-labeled
2-dimensional subsets of parameters defined by $p^m=q^{c+2}$. They are
thus
generalized $\qpvir$ algebras at $c=-2+m\frac{\ln p}{\ln q}$.
In such a frame the closed algebraic relation (\ref{eq51}) may
acquire a crucial importance as a $q$-deformation of the Virasoro-current
commutations relations.
A better understanding of the undeformed limit $q \rightarrow 1$ would
help
us to clarify this interpretation if one could indeed identify the
standard
Virasoro--Kac-Moody structure in such a limit. The difficulty lies in the
correct definition of this limit for the generators $L^\pm(z)$ and $t(z)$
which
should be consistent with such an interpretation.
As for the previously mentioned central-extension
problem, a help could come from an explicit bosonization of the
elliptic algebra, as was done for $U_q(\widehat {sl}(N)_c)$ in
\cite{AOS}.

\sect {Conclusion}

We have studied some aspects of the elliptic quantum algebra $\elp$,
in order to show its importance as a generalization of the quantum affine
algebra $U_q(\widehat {sl}(2)_c)$. We have seen that the introduction of
$\elp$ permits to incorporate in the context of a deformed affine algebra
the
$q$-deformed Virasoro algebra introduced by \cite{SKAO}, which is the
symmetry of trigonometric Ruijsenaars--Schneider model. This is obtained
on the particular surface of the space of parameters of $\elp$ given by
the equation $p=q^{c+2}$.
We expect that the other $\qpvir$ algebras constructed on the surfaces
$p^m=q^{c+2}$, $m\not =1$, will provide us with the mathematical
structure required to study other relativistic integrable models,
characterizing elliptic quantum algebras as a general framework for the
description of symmetries in (quantum) relativistic mechanics.

\end{document}